\documentclass{birkjour}

\usepackage{amsmath,amsthm,amssymb,amscd,amsfonts}
\usepackage{mathtools,slashed}
\usepackage{amssymb,latexsym}
\usepackage{tikz}
\usetikzlibrary{matrix,arrows}
\usepackage{graphicx}
\usepackage[hyperindex]{hyperref}
 \newtheorem{thm}{Theorem}[section]

 \theoremstyle{definition}
 \newtheorem{defn}[thm]{Definition}
 \theoremstyle{remark}

 \numberwithin{equation}{section}

\begin{document}

%
%
%
%
%
%
%
%
%

\title[Equivariant homologies for operator algebras]
 {Equivariant homologies for operator \\
 	algebras}

\author[M. Amini]{Massoud Amini}

\address{%
Department of Mathematics, Tarbiat Modares University\\
 P.O. Box 14115-175\\
Tehran\\
Iran}

\email{mamini@modares.ac.ir}

\author[A. Shirinkalam]{Ahmad Shirinkalam}
\address{%
	Department of Mathematics, Tarbiat Modares University\\
	P.O. Box 14115-175\\
	Tehran\\
	Iran}
\email{shirinkalam\_a@aut.ac.ir}



\begin{abstract}
This is a survey of a variety of equivariant (co)homology theories for operator algebras. We briefly discuss a background on equivariant theories, such  as equivariant $K$-theory and equivariant cyclic homology. As the main focus, we discuss a notion of equivariant $ L^2 $-cohomology and equivariant $ L^2 $-Betti numbers for  subalgebras of a von Neumann algebra. For graded  $C^*$-algebras  (with grading over a group) we elaborate on a notion of graded $ L^2 $-cohomology and  its relation to  equivariant $L^2$-cohomology.
\end{abstract}

\maketitle
\section{Introduction}
Homology theories are generally regarded as a method of constructing algebraic invariants of other structures, with a wide range of applications in geometry, analysis  and algebra. 

The first historical instance of a homology theoretical consideration is probably the Euler polyhedron formula (now called Euler characteristics), rigorously proved  by Cauchy in 1811 (for any convex polyhedron). It was  Riemann who first 's defined the genus (and connectedness numerical invariants) in 1857 and Betti who showed   the independence of homology numbers from the choice of basis in 1871.

Henri Poincar\`{e}  invented
 the fundamental group and initiated the study of the so called algebraic varieties and their homology theory. He was then led to Betti numbers through solving certain vector field problems  (to determine the number of holes in the
configuration space). One could make a simple observation of how  holes (and so betti numbers) are  significant by examining the behavior of solutions of differential equations, say  on the torus compared to the sphere. Poincar\`{e} published his ideas of analysis situs in 1895, where he introduced the notion of homology classes. To classify  possible configurations of orientable cycles, he used Betti numbers (as refinements of the Euler characteristics). He also noticed that classifying non-orientable cycles requires torsion coefficients.

The history goes on with great names, such as Alexander (with theory of knots) Veblen,  Brouwer, as well as van Kampen and Lefschetz.  It was Emmy Noether who first recognized that the Betti
numbers and torsion coefficients are numerical invariants of isomorphism
classes of finitely-generated abelian groups.

Algebraic counterparts of the theory came later. As one of the first instances of such algebraic homology theories, the Hochschild (co)homology was introduced by Gerhard Hochschild  in 1945 for algebras over a field, and extended to algebras over more general rings by Henri Cartan and Samuel Eilenberg  around 1956.

Cyclic (co)homology and cyclic cohomology which generalize the de Rham (co)homology of manifolds were independently introduced by Boris Tsygan (homology)  and Alain Connes (cohomology)  in the 80's. These are naturally related to the de Rham theory, Hochschild (co)homology and K-theory. The later development are, among  others, due to Max Karoubi, Mariusz Wodzicki and Jean-Louis Loday.

The origin of the equivariant homology theory (motivated by group extensions) goes back to Whitehead in 1950. In  algebraic topology, this provides an invariant for classification of topological spaces with a (given) group action. This is motivated by  gauge theory in physics, which  is a  field theory where the Lagrangian is invariant under certain Lie groups of local transformations. It is shown that all homology theories have very natural and successful equivariant counterparts.

The main focus of the present survey is a notion of  equivariant $ L^2 $-cohomology for operator algebras. This is discussed in more details in Section 3. In this section we briefly review some of the other existing equivariant theories for operator algebras. These in particular include  equivariant $K$-theory and equivariant cyclic homology. We also discuss equivariant Hochschild cohomology in somewhat more technical details and return to it again in the last section, where we discuss equivariant $ L^2 $-cohomology and equivariant $ L^2 $-Betti numbers.

\subsection{Topological Background} The main idea in singular cohomology
is to give a functor from the category of
topological spaces and continuous maps to the
category of graded rings and homomorphisms.
When there is also  a group
action on the topological space (with appropriate continuity properties) one may ask about the existence of a functor which encodes
both the topology  and the algebra (the action). Equivariant cohomology is supposed to do this.

The origin of equivariant cohomology goes back to  1959 when Borel defined equivariant
singular cohomology using the so called Borel construction \cite{tu}.
In 1950, Cartan studied the action of a compact connected Lie group on a manifold
and constructed a differential complex
out of the differential forms on the manifold and the Lie
algebra of the acting group (which gives the real
equivariant singular cohomology in modern terms). When the group is trivial, this reduces to the de Rham complex
of smooth differential forms which computes the real
singular cohomology of the manifold. The equivariant cohomology proved later to be a strong computational tool, when it was used
in
the equivariant integration formula of Atiyah-Bott
(1984) and Berline-Vergne (1982) \cite{tu}. In turn, the equivariant integration formula is used to show that the
stationary phase approximation formula is exact
for a symplectic action (Atiyah-Bott 1984), to calculate
the number of rational curves in a quintic
threefold (Kontsevich 1995, Ellingsrud-Strømme
1996), and to calculate the characteristic numbers of
a compact homogeneous space (Tu 2010) \cite{tu}. Following Tu, we refer the reader to the expository articles by Bott \cite{b} and Vergne \cite{v} for further applications.

\subsection{Equivariant $K$-Theory} Equivariant K-theory is introduced by Atiyah in 60's. The basic idea is to replace vector bundles with $G$-vector bundles over a $G$-space $X$. The Bott periodicity \cite{ab} (in the complex setting) is generalized by Atiyah and Segal, and various applications are put forward \cite{s}. Among the most sophisticated applications, one should mention the Atiyah-Singer index theorem (c.f. \cite{a} and references therein).

Also group actions on varieties is studied by Goresky, Kottwitz, and MacPherson \cite{gkm} to give an approach to equivariant
cohomology which  uses the combinatorial data of
the orbits of the  action (moment map theory) to  interpret
topological data of the attracting sets for fixed points \cite{ty}. This is successfully  extended to other cohomology theories, including equivariant intersection
homology \cite{bm}, equivariant K-theory \cite{rk}, and  is used to calculate the equivariant
cohomology rings of Grassmannians \cite{kt} and regular varieties \cite{bc}.

\subsection{Equivariant Cyclic Homology}
Cyclic homology was introduced by Connes in 80's as a noncommutative analogue
of de Rham homology \cite{c1}. On the other hand, cyclic cohomology happens to be a   natural extension of the Chern
character in K-homology \cite{c2}. Loday and Quillen  computed the cyclic homology
of the algebra of regular functions on smooth varieties \cite{lq} and the same is done by Burghelea for group rings \cite{bu}.

In noncommutative geometry, the equivariant theory is basically formulated through crossed product algebras. The cyclic homology for crossed-product algebras
is known in the case of actions of finite groups (c.f. \cite{p} and references therein). For general group
actions, Feigin and Tsygan constructed a spectral sequence converging to cyclic homology \cite{ft}. For group actions on manifolds, the case of
discrete proper actions was settled by Baum and Connes \cite{bc1}. For general group actions,
the periodic cyclic homology of the homogenous component is computed in terms of equivariant
homology via the construction of an explicit cochain map (see \cite{p}, for this and other advances in this direction).

\subsection{Equivariant Hochschild Cohomology} The equivariant Hochschild cohomology theory for Banach and operator algebras is developed by K. K. Jensen \cite{j}. In many mathematical and physical applications, the structure involved carries a naural group action. The most natural examples come from  symmetries or  dynamical evolutions in a physical context. The main objective of the equivariant Hochschild cohomology theory is to describe and classify the structures which contain such group actions, and to give an invariant for associated geometric constructs.

In \cite{j}, Jensen introduces the cohomology first in a purely algebraic setting and then adapt it to the context of Banach and operator algebras. He also calculates some low-dimensional equivariant Hochschild cohomologies. For example, he showes that If $ \mathcal{A} $ is a unital abelian $C^*$-algebra and $ \alpha $ is an action of a group $ G $ on $ \mathcal{A} $, then $ H^2_{G}(\mathcal{A},\mathcal{A})=0 $ \cite[Corollary II.1.10]{j}. He also verifies the reduction dimension formula \cite[Corollary II.2.6]{j} as
$H^{n+p}_G(\mathcal{A}, X)\simeq H^{n}_G(\mathcal{A}, BL^p(\mathcal{A}, X)). $

The relation between ordinary and equivariant cohomology is also investigated. Let $ \mathcal{A} $ be a Banach algebra and
$ G $ be a locally compact group acting on $ \mathcal{A} $. Then the map $ H^{n}_G(\mathcal{A}, X)\rightarrow H^{n}(\mathcal{A}, X) $ is injective, whenever $ G $ is compact (amenable) and $ X $ is a Banach (dual) module \cite[Proposition II.3.2]{j}. Also if  the action $ \alpha $ of $ G $ on $ \mathcal{A} $ is inner and the Banach algebra generated by $ \{\alpha_g | g \in G\} $ is amenable and the module actions on $ X $ have norm at most 1, then the map $ H^{n}_G(\mathcal{A}, X)\rightarrow H^{n}(\mathcal{A}, X) $ is also surjective.
In particular, if $ G $ is an amenable group that acts innerly on a von Neumann algebra $ \mathcal{M} $, then $ H^{n}_G(\mathcal{M}, \mathcal{M})$ and $ H^{n}(\mathcal{M}, \mathcal{M}) $ are isomorphic \cite[Lemma II.3.10]{j}.

Jensen also showes that if $\mathcal{A}  $ is a UHF $C^*$-algebra and $ G $ has product action on $ \mathcal{A} $, then $H^{n}_G(\mathcal{A}, X)=0  $ for every dual module $ X $ \cite[Corollary III.2.8]{j}. For a  II$_1 $-factor $ \mathcal{M} $ with a discrete group action,
$ \mathcal{M} $ is injective whenever $ H^{1}_G(\mathcal{M},X)=0 $ for every dual module $ X $ \cite[Theorem III.3.7]{j}. We return to this theory in Section 3.1, where we give more details needed for equivariant $L^2$-cohomology.

\section{$ L^2 $-Cohomology}

\subsection{Basic Notions}
Throughout the rest of the paper, for an algebra $ \mathcal{A} $
we denote its opposite algebra  by $ \mathcal{A}^\circ $ and the enveloping algebra associated with $ \mathcal{A} $ by $ \mathcal{A}^e=\mathcal{A}\otimes \mathcal{A}^\circ $.

Also,
$ \mathcal{A}_1 $ denotes the closed unit ball of a $C^*$-algebra $ \mathcal{A}$.
The weak (respectively, strong) operator topology is denoted by WOT (respectively, SOT) and UWOT denotes the ultra-weak operator topology on $ B(\mathcal{H}) $, the space of bounded linear operators on a Hilbert space $ \mathcal{H} $.

Let $X$ be a Banach $\mathcal{A}$-bimodule. Then, $ X $ is said to be dual if it is a dual of a Banach space and for each $ a\in \mathcal{A} $, the  maps $X \rightarrow X; x\mapsto a\cdot x  $ and $ x\mapsto x\cdot a  $
are weak* continuous. If in addition, $ \mathcal{A} $ admits a weak* topology (for example whenever $ \mathcal{A} $ is a von Neumann algebra), and for every $ x\in X $ the  maps $\mathcal{A}  \rightarrow X; a\mapsto a\cdot x  $ and $ a\mapsto x\cdot a  $ are weak* continuous, then $ X $ is called a normal dual module.

Let $X$ be a Banach $\mathcal{A}$-bimodule. Let $ BL^0(\mathcal{A},X)=X $ and for each $ n \in \mathbb{N} $,  denote by $ BL^n(\mathcal{A},X) $ the space of all bounded $ n $-linear maps from $ \mathcal{A}^n $ into $ X. $

For a $C^*$-algebra $\mathcal{A}$, a   map $ \tau:\mathcal{A}^+\rightarrow [0 , \infty] $   is called
tracial (or a trace) if $ \tau(ab)=\tau(ba) $, $ \tau(a+b)=\tau(a)+\tau(b) $ and $ \tau(\lambda a)=\lambda \tau(a), $ for all $a,b \in \mathcal{A}^+$ and $ \lambda\geq 0 $.
A trace  $\tau$
is said to be faithful if $a=0$ whenever $\tau(a^*a)=0$, for  $a\in \mathcal{A}$.
Each faithful trace on $\mathcal{A}$ induces a norm $ \mathcal{A}rrowvert . \mathcal{A}rrowvert_{\tau} $ on $\mathcal{A}$  defined by $\mathcal{A}rrowvert a \mathcal{A}rrowvert^2_{\tau}=\tau(a^*a) \ (a\in \mathcal{A})$.
A trace  $\tau$  on a von Neumann algebra $\mathcal{M}  $ is said to be normal if $ a_i \uparrow a $ in the SOT of $\mathcal{M}$, then $ \tau (a_i) \uparrow \tau (a)$.

Let $ \mathcal{M}\subseteq B(\mathcal{H}) $ be a von Neumann algebra. A closed densely defined operator $ T:Dom(T)\subseteq\mathcal{H}\rightarrow \mathcal{H} $ is said to be affiliated with $ \mathcal{M} $ if $ Tu=uT $ for all unitary elements in the commutant $ \mathcal{M}^\prime. $ The set of
all  operators affiliated with $ \mathcal{M} $ is denoted by $ \mathcal{U}(\mathcal{M})$.
If $ T $ is a closed densely defined operator on $\mathcal{H}  $ with the polar decomposition $ T=v\arrowvert T\arrowvert $, then $ T\in \mathcal{U}(\mathcal{M}) $ if and only if $ v\in \mathcal{M} $ and $ \arrowvert T\arrowvert = \int_0^\infty \lambda \, dE_\lambda \in \mathcal{U}(\mathcal{M}) $. In this case, the spectral projections $ E_\lambda $ belong to $ \mathcal{M}$.

For a finite von Neumann algebra $ \mathcal{M}$, the set $ \mathcal{U}(\mathcal{M}) $ is an involutive algebra under the natural operations. The interested reader will find more details in \cite{be} and \cite{d1}.

E. Nelson in \cite{n} showed that if  $ \mathcal{M}\subseteq B(\mathcal{H}) $ is a tracial von Neumann algebra with a faithful, normal, finite trace $ \tau $, then
$ \mathcal{U}(\mathcal{M}) $  carries a translation-invariant topology  which is called the measure topology and is generated by the family $ \{ N(\epsilon,\delta)\}_{\epsilon,\delta >0} $ of neighborhoods of 0  of the form
$$ N(\epsilon,\delta)=\{  a\in \mathcal{M}; \hbox{for some projection} \,p\in  \mathcal{M},\| ap \| \leq\epsilon \, \hbox{and}\: \tau (1-p)\leq\delta \}.$$
He also showed   that one can identify
$ \mathcal{U}(\mathcal{M}) $ with the closure of $\mathcal{M}  $ in this topology \cite[Theorem 4]{n}. More precisely, each $ T \in \mathcal{U}(\mathcal{M}) $ with the spectral decomposition $ v\int_0^\infty \lambda \, dE_\lambda $ is the limit of a sequence $ T_n=v\int_0^n \lambda \, dE_\lambda \in \mathcal{M} $ in the measure topology.

We use this  fundamental fact in Section 3 to extend an action of a group $ G $
on $ \mathcal{M} $ to an action of $ G $
on $ \mathcal{U}(\mathcal{M})$.

Let $ \mathcal{A} $ be a unital algebra with a unit 1.  Then, $\mathcal{A}^e  $ is an $ \mathcal{A} $-bimodule with the actions
\begin{equation}\label{a1}
a\cdot (b\otimes c^\circ)=ab \otimes c^\circ \quad \hbox{and} \quad (b\otimes c^\circ)\cdot a = b \otimes (ca)^\circ \qquad (a,b,c \in \mathcal{A}).
\end{equation}
If $ \mathcal{A} $ is a Banach algebra, these actions extend by linearity and continuity.

In general, there is a one to one correspondence between bimodules over $ \mathcal{A} $ and left (right) modules over $ \mathcal{A}^e $. Indeed, if $ X $ is an $ \mathcal{A} $-bimodule, then it is a left $ \mathcal{A}^e $-module via $ (a\otimes b^\circ)\cdot x=a\cdot x \cdot b  $ (or it is a right $ \mathcal{A}^e $-module via $ x \cdot (a\otimes b^\circ)=b\cdot x \cdot a  $) for  $ a,b \in \mathcal{A} $ and $ x \in X $. Conversely, if
$ X $ is a left $ \mathcal{A}^e $-module, then it is an $ \mathcal{A} $-bimodule with the actions $ a\cdot x = (a \otimes 1^\circ) \cdot x $  and $ x \cdot b=(1\otimes b^\circ)\cdot x $ for every $ a,b \in \mathcal{A} $ and $ x \in X $.
The right module actions  are defined similarly.

For a von Neumann algebra $\mathcal{M} $, its enveloping von Neumann algebra
$ \mathcal{M}\bar{\otimes}\mathcal{M}^\circ $ is denoted by $ \mathcal{M}^e $. Since $\mathcal{U}(\mathcal{M}^e)  $ is a left $ \mathcal{M}^e $-module with the multiplication, by the above argument, $\mathcal{U}(\mathcal{M}^e)  $ is a $ \mathcal{M} $-bimodule with the module actions
\begin{equation}\label{a2}
m \cdot T=(m\otimes 1^\circ) T, \quad T\cdot n = (1\otimes n^\circ) T \qquad(m,n \in \mathcal{M}, \:T \in \mathcal{U}(\mathcal{M}^e)).
\end{equation}

\subsection{$ L^2 $-Cohomology for von Neumann Algebras}

Finding a suitable homology in the context of von Neumann algebras goes back to the works of Johnson, Kadison and Ringrose (1971-72) and then to Sinclair and Smith (1995). Their goal was to find a powerful invariant to distinguish von Neumann algebras.

The theory of $ L^2 $-homology for von Neumann algebras is introduced by A. Connes and D. Shlyakhtenko \cite{cs} following the works of D. Gaboriau (2002) in the field of ergodic equivalence relations, where he showed that all the $ L^2 $-Betti numbers are the same for those discrete groups that can generate isomorphic ergodic measure-preserving equivalence relations. Connes and Shlyakhtenko used the theory of correspondences together with the algebraic description of $ L^2 $-Betti numbers given by L\"{u}ck, to define the $ k $-th  $ L^2 $-homology of a von Neumann algebra
$ \mathcal{M} $ by $$ H^{(2)}_k(\mathcal{M})=H_k(\mathcal{M},\mathcal{M}\bar{\otimes}\mathcal{M}^\circ), $$
where $ H_k $ means the algebraic Hochschild homology. Then the corresponding  $ L^2 $-Betti numbers of $ \mathcal{M} $ are defined to be $ \beta^{(2)}_k( \mathcal{M})=\dim_{\mathcal{M}\bar{\otimes}\mathcal{M}^\circ}H^{(2)}_k(\mathcal{M})$.
The dimension function used here is the extended dimension function of L\"{u}ck, that is explained in more details in the next section.

Connes and Shlyakhtenko investigate the relation between  $ L^2 $-Betti numbers of a discrete group $ \Gamma $ with $ L^2 $-Betti numbers of its group algebra and  showed that $\beta^{(2)}_k(\Gamma)=\beta^{(2)}_k(\mathbb{C}\Gamma, \tau),  $ where $ \tau $ is the standard trace on the group algebra $ \mathbb{C}\Gamma $ \cite[Proposition 2.3]{cs}. For a  II$_1 $-factor $ \mathcal{M} $ they computed its
$ L^2 $-Betti numbers. They showed that the zeroth Betti number is 0 and  each $ L^2 $-Betti number of $ \mathcal{M} $ is a limit of $ L^2 $-Betti numbers of its sub-complexes \cite[Lemma 2.2]{cs}.
They also showed that for a von Neumann algebra $ (\mathcal{M},\tau) $ with a diffuse
center, $ \beta^{(2)}_1(\mathcal{M},\tau)=0 $ \cite[Corollary 3.5]{cs}. In particular, if $ \mathcal{M} $ is commutative, then $ \beta^{(2)}_k(\mathcal{M})=0, $ for every $ k\geq1 $ \cite[Corollary 5.4]{cs} (the same result is obtained by A. Thom for all von Neumann algebras with a diffuse
center \cite[Theorem 2.2]{t}). Finally they give connections to free probability theory including the inequality between the microstates and  microstates-free entropy.

In Section 3 we show how to combine the Jensen equivariant theory and Connes-Shlyakhtenko $L^2$-theory to give a numerical invariant
in the context of dynamical systems of operator algebras.
To do this, we first introduce a notion of equivariant $ L^2 $-cohomology and equivariant $ L^2 $-Betti numbers for a subalgebra of a von Neumann algebra. Then, considering a grading of a  group $ G $ on a $C^*$-algebra $ \mathcal{A} $, we give a definition of the graded $ L^2 $-cohomology and the associated  $ L^2 $-Betti numbers for $ \mathcal{A}$.
We investigate the relation between graded $L^2$-cohomology of $ \mathcal{A} $ and the equivariant $L^2$-cohomology of $ \mathcal{A} $, when the group $ G $ is abelian.

\subsection{The Dimension Function}
Here we give a brief review of the generalized Murray-von Neumann dimension function. The interested reader may find more details in \cite[Chapter 6]{lu}.

Let $ R $ be a unital ring and $ X  $ be an $ R $-module. The dual module $ X^* $  of $ X $ is  Hom$_R(X,R), $ where the $ R $-multiplications are given by $(fr)(x)=f(x)r  $
and $ (rf)(x)=rf(x), $ for $ f\in X^*, x\in X $ and $ r\in R$.
If $ Y $ is an $ R $-submodule of $ X$, then the closure of $ Y $ in $ X $ is defined to be
$$ \overline{Y}=\{x\in X | f(x)=0\: \hbox{for all}\: f\in X^* \;\hbox{with} \;Y\subseteq \ker (f) \}. $$

For an $ R $-module $ X, $ define the  submodule $ \mathrm{T}X $ and the quotient module $\mathrm{P}X $ by
$$ \mathrm{T}X:=\{x\in X| f(x)=0\: \hbox{for all}\: f\in X^* \}=\overline{\{0\}} $$
and
$$\mathrm{P}X:=X/\mathrm{T}X.  $$
Note that $\mathrm{T}\mathrm{P}X=0, \mathrm{P}\mathrm{P}X=\mathrm{P}X, X^*=(\mathrm{P}X)^*,$
and that $ \mathrm{P}X=0 $ is equivalent to $ X^*=0. $ If $ X $ is a finitely generated (f.g.) $ R $-module, then $ \mathrm{P}X $ is finitely generated projective and $ X\cong\mathrm{T}X\oplus \mathrm{P}X $.

Let $ \mathcal{M}$ be a finite von Neumann algebra with a fixed normal, faithful state $ \tau $. The tracial functional $ \tau_n $ on $ M_n(\mathcal{M}) $ given by $\tau_n(\{ a_{ij}\}_{i,j=1}^{n})=\sum_{i,j=1}^{n}\tau(a_{ii})  $ is positive and faithful (but not a state, since it takes the value $ n $ on the unit matrix).
If $ P $ is a finitely generated projective $ \mathcal{M} $-bimodule, then it is isomorphic to $ \mathcal{M}^nA $ for some idempotent matrix $ A\in M_n(\mathcal{M}) $.
Now the Murray-von Neumann dimension of $ P $ is defined to be
$$\dim_{\mathcal{M}}(P)=\tau_n(A)\in [0,\infty).  $$
This definition is independent of the choice of the matrix $ A $.
If $ X $ is an arbitrary $ \mathcal{M} $-bimodule, then the dimension of $ X $ is defined by
$$\dim^\prime_{\mathcal{M}}(X)=\sup \{\dim_{\mathcal{M}}(P)|P \hbox{ is a f.g. projective submodule of }X \}\in [0,\infty].  $$
By \cite[Theorem 6.7]{lu}$, \dim^\prime $ is the only dimension function that extends
the Murray-von Neumann dimension and satisfies the following properties;

(1) (Extension Property) If $ P $ is a finitely generated projective $ \mathcal{M} $-module, then $$\dim^\prime_{\mathcal{M}}(P)=\dim_{\mathcal{M}}(P);  $$

(2) (Additivity) If $ 0\rightarrow X_0 \rightarrow X_1 \rightarrow X_2\rightarrow 0 $ is an exact sequence of modules, then $$ \dim^\prime_{\mathcal{M}}(X_1) = \dim^\prime_{\mathcal{M}}(X_0) + \dim^\prime_{\mathcal{M}}(X_2);  $$

(3) (Continuity) If $ Y $ is a submodule of the finitely generated module $ X $, then
$$ \dim^\prime_{\mathcal{M}}(Y)=\dim^\prime_{\mathcal{M}}(\overline{Y}); $$

(4) If $ X $ is finitely generated $ \mathcal{M} $-module, the
$$\dim^\prime_{\mathcal{M}}(X)=\dim_{\mathcal{M}}(\mathrm{P}X) \; \hbox{and} \: \dim^\prime_{\mathcal{M}}(\mathrm{T}X)=0. $$

\section{Equivariant $ L^2 $-Cohomology}
\subsection{Equivariant Hochschild Cohomology}

Following  \cite{j}, we start with the notion of  action of a  locally compact group on a $ C^* $-algebra.

\begin{defn} \cite[Definition I.1.1]{j}
	A  group $G$ is said to act on an algebra $\mathcal{A}$ if there is
	a homomorphism $ \alpha$ from $G$ into  the group of automorphisms of $\mathcal{A}$, that is,
	a map $ \alpha : G\rightarrow Aut(\mathcal{A});g \mapsto \alpha_g $ with the following properties:
	
	(i) for each $ g \in G $ the map $\alpha_g  $ is  bijective and linear,
	
	(ii) $ \alpha_g(ab)=\alpha_g(a) \alpha_g(b) \ (a,b \in \mathcal{A}), $

	(iii) $ \alpha_{gh}=\alpha_g \alpha_h \ (g,h \in G). $
\end{defn}
We denote this action by $ G \curvearrowright _{\alpha} \mathcal{A} $. Also, the fixed point algebra of $ \mathcal{A} $ with respect to $ \alpha $ is the set $ \mathcal{A}^\alpha =\{a \in A | \alpha_g(a)=a \ (g \in G)\} $.

\begin{defn} \cite[Definition I.1.2]{j}\label{equi}
	Let $ G \curvearrowright _{\alpha}\mathcal{A} $. Suppose $ X  $ is an  $\mathcal{A}$-bimodule.
	Then, $ X $ is said to be an equivariant $ (\alpha,\beta) $-$ G$-$\mathcal{A}$-module if
	there is a map
	$ \beta : G\rightarrow Lin(X);g \mapsto \beta_g $ with the following properties:
	
	(i) for each $ g \in G $ the map $\beta_g  $ is a  linear bijective,
	
	(ii) $ \beta_g(a\cdot x)=\alpha_g(a) \cdot \beta_g(x) \ \hbox{and}  \:\beta_g(x\cdot a)= \beta_g(x) \cdot  \alpha_g(a) \ (a \in \mathcal{A}, x \in X), $

	(iii) $ \beta_{gh}=\beta_g \beta_h \ (g,h \in G). $
\end{defn}
We denote this action by $ G \curvearrowright _{\beta}X $.

If $ G \curvearrowright _{\alpha}\mathcal{A} $, then $ \mathcal{A} $ provided with the module actions of left and right multiplication with $ \beta_g = \alpha_g $ is a
Banach equivariant $ G$-$\mathcal{A}$-module.

On the other hand, let $\mathcal{M}$ be a  von Neumann algebra and let $ G \curvearrowright _{\alpha}\mathcal{M} $. By the module structure of $\mathcal{M}^e  $ described in (\ref{a1}) and by an argument similar to \cite[I.1.4 (a)]{j}, $ G $ acts  on $ \mathcal{M}^e $ via the map
$ \tilde{\alpha}_g:\mathcal{M}^e\rightarrow \mathcal{M}^e$ defined by $\tilde{\alpha}_g(a\otimes b^\circ)= \alpha_g(a)\otimes \alpha_g^\circ(b) \quad (a,b \in \mathcal{M}, g \in G )$, extended by linearity and continuity. Hence $ \mathcal{M}^e $ is an equivariant
$ (\alpha,\tilde{\alpha})$-$G$-$\mathcal{M}$-module.

Let $\mathcal{M}$ be a  von Neumann algebra with a finite, faithful normal trace $ \tau $ and let $ G \curvearrowright _{\alpha}\mathcal{M} $.
This trace $ \tau $ is called invariant under this action if for each $ g \in G $ and $ a \in \mathcal{M} $ we have $ \tau(\alpha_g(a))=\tau(a) $.

Let us observe that if
$\mathcal{M}$ is a  von Neumann algebra with a finite, faithful,  normal trace $ \tau $ and  $ G \curvearrowright _{\alpha}\mathcal{M} $ with the extra property that
$\tau  $ is invariant under the action $ \alpha $, then there is an action $ \beta:G \rightarrow Lin(\mathcal{U}(\mathcal{M}^e)) $ which makes $ \mathcal{U}(\mathcal{M}^e) $ an equivariant $ (\alpha,\beta)$-$G$-$\mathcal{M}$-module.

To see this,	by the actions described in (\ref{a1}), $ \mathcal{M}^e $ is a Banach $ \mathcal{M} $-bimodule and
 as above, we have an action $ \tilde{\alpha}_g:\mathcal{M}^e\rightarrow \mathcal{M}^e \quad (g\in G), $ making $ \mathcal{M}^e $ into an equivariant $ \mathcal{M} $-module.
Since $ \mathcal{M}^e $ is dense in $ \mathcal{U}(\mathcal{M}^e) $ in the measure topology, for every $ T\in \mathcal{U}(\mathcal{M}^e) $ there is a sequence $ (T_n)\subseteq \mathcal{M}^e $ such that $ T_n\rightarrow T $
in the measure topology. Now we define
an action
$ \beta:G \rightarrow Lin(\mathcal{U}(\mathcal{M}^e)) $  by
$ \beta_g(T)=\lim_n \tilde{\alpha}_g(T_n) \ (T\in \mathcal{U}(\mathcal{M}^e), g \in G )$, where the limit is  taken in the measure topology. If $ T=0 $, then there is a sequence
$ (T_n)\subseteq \mathcal{M}^e $ such that $ T_n\rightarrow 0 $. Hence by the definition of 0 neighborhoods in the measure topology, we have
$$\forall \epsilon , \delta >0 \quad   \exists n_0, \quad T_n\in N(\epsilon,\delta) \quad\hbox{whenever} \quad n>n_0.$$
Since $ \tau  $ is invariant under $ \alpha $, the
finite, faithful, normal trace
$ \tau\bar{\otimes}\tau^\circ $ defined by $ \tau\bar{\otimes}\tau^\circ(a\otimes b^\circ) =\tau (a)\tau^\circ(b^\circ)$
is invariant under $ \tilde{\alpha} $. This  shows that
$ T_n\in N(\epsilon,\delta) $ if and only if $ \tilde{\alpha}_g(T_n)\in N(\epsilon,\delta)$, so that $ \tilde{\alpha}_g(T_n)\rightarrow 0 $. Thus
$ \beta_g(T)=0 $. This argument justifies that  the action $ \beta $ is well-defined.

For every $ a, b \in \mathcal{M} $ by the actions of
$ \mathcal{M} $ on $ \mathcal{U}(\mathcal{M}^e) $ described in (\ref{a2})
we have
\begin{align*}
\beta_g(a\cdot T)=\beta_g((a\otimes 1^\circ) T)
&=  \lim_n \tilde{\alpha}_g((a\otimes 1^\circ) T_n)   \\
&=  \lim_n (\alpha_g(a)\otimes 1^\circ) \tilde{\alpha}_g(T_n) \qquad (\tilde{\alpha} \;\hbox{is equivariant}) \\
&=  (\alpha_g(a)\otimes 1^\circ)\cdot\lim_n  \tilde{\alpha}_g(T_n) \qquad  \\
&=\alpha_g(a)\cdot\beta_g(T).
\end{align*}
A similar argument shows that $ \beta_g(T\cdot b)=\beta_g(T)\cdot \alpha_g(b).$

If $ g, h\in G  $, then
\begin{align*}
\beta_{gh}(T)&= \lim_n \tilde{\alpha}_{gh}(T_n)=\lim_n \tilde{\alpha}_{g}(T_n)\tilde{\alpha}_{h}(T_n)\\&=\lim_n \tilde{\alpha}_{g}(T_n)\lim_n \tilde{\alpha}_{h}(T_n)=\beta_{g}(T)\beta_{h}(T),\end{align*}
thus $ \beta $ is a homomorphism.

Finally, we observe that $ \beta_{g} $ is bijective.
For each $ g\in G $ and $ T \in \mathcal{U}(\mathcal{M}^e) $, $ \beta_{g}(T) $ has a dense domain $ D \subseteq \mathcal{L}^2(\mathcal{M}\bar{\otimes}\mathcal{M}^\circ) $.
We show that for such $ g $, $ \beta_{g^{-1}}(T) $ with domain $ D^\prime $ is the inverse of
$ \beta_{g}(T) $, that is, $ \beta_{g^{-1}}(T)\beta_{g}(T)\xi=T\xi, $ for all  $\xi \in D \cap Dom(T)   $ with $ \beta_{g}(T)\xi \in D^\prime $
and $ \beta_{g}(T)\beta_{g^{-1}}(T)\eta=T\eta, $ for all $ \eta \in D^\prime \cap Dom(T) $ with $\beta_{g^{-1}}(T)\eta \in D  $. If $ (T_n)_n $ is a sequence in $ \mathcal{M}^e $ converging to $ T $ in the measure topology, then
$$ \beta_{g}(T)\beta_{g^{-1}}(T)=\lim_n \tilde{\alpha}_{g}(T_n)\lim_n \tilde{\alpha}_{g^{-1}}(T_n)=\lim_n \tilde{\alpha}_{gg^{-1}}(T_n)=\lim_n (T_n)=T, $$
and similarly, $ \beta_{g^{-1}}\beta_{g}(T)= T $. Hence $ \mathcal{U}(\mathcal{M}^e) $ is an equivariant $ (\alpha,\beta)$-$G$-$\mathcal{M}$-module.

Next let us recall the definition of an equivariant cochain.

\begin{defn}\cite[Definition I.1.5, I.1.7]{j}
	Let $ G, \mathcal{A}, X, \alpha $ and $ \beta $ be as in Definition \ref{equi}.
	For each $ n\geq 1 $ define
	$$C^n_G(\mathcal{A},X)=\{ T\in L^n(A,X)| T(\alpha_g(a_1),...,\alpha_g(a_n))=  \beta_g(T(a_1,...,a_n)) \ (g\in G) \}.  $$
	The elements of $ C^n_G(\mathcal{A},X) $ are  called the $ (\alpha,\beta) $-equivariant $ n $-cochains (of $ \mathcal{A} $ with
	coefficients in $ X $).
	
	The equivariant zero-cochains are defined separately as follows
	$$C^0_G(\mathcal{A},X)=\{ x\in X|  \beta_g(x)=x\ (g\in G) \}.  $$
	Note that  the usual coboundary operator $ \delta^n:C^n(\mathcal{A},X)\rightarrow C^{n+1}(\mathcal{A},X) $
	maps $ C^n_G(\mathcal{A},X) $ into $ C^{n+1}_G(\mathcal{A},X) $. We denote the restriction map by $ \delta^n_G $.
	
	The equivariant cocycle and coboundary
	groups are defined, respectively, by
	$$Z^n_G(\mathcal{A},X)=\ker \delta^n_G, \quad  B^n_G(\mathcal{A},X)=\hbox{ran}\: \delta^{n-1}_G$$
	and the equivariant cohomology groups are defined as the quotients
	$$ H^n_G(\mathcal{A},X)=Z^n_G(\mathcal{A},X)/B^n_G(\mathcal{A},X). $$
\end{defn}

\subsection{Equivariant $ L^2 $-Cohomology}	
Next we define the notion of equivariant $ L^2 $-cohomology
for a tracial $ * $-algebra.

Let $ (\mathcal{A},\tau) $ be a tracial $ * $-algebra satisfying in the following two conditions:

1. $ \tau(a^*a)\geq 0 $ for all $ a \in \mathcal{A} $; and

2. $ \forall b \in \mathcal{A} \: \exists C>0: \ \tau(a^*b^*ba)\leq C\tau(a^*a) \quad (a\in \mathcal{A}).$

Suppose that  $ \mathcal{M}=W^*(\mathcal{A})\subseteq B(\mathcal{L}^2(\mathcal{A},\tau)) $
be its enveloping von Neumann algebra. We adapt the following definition due to A. Thom \cite[Definition 3.6]{t}.

\begin{defn}
	With the above assumptions, the $ k $-th $ L^2 $-cohomology of $ \mathcal{A} $ is defined by $$ H^{(2),k}(\mathcal{A})=H^k(\mathcal{A}, \mathcal{U}(\mathcal{M}^e)),  $$
	where $H^k(\mathcal{A}, \mathcal{U}(\mathcal{M}^e))  $ is the $ k $-th  Hochschild cohomology of $ \mathcal{A} $ with coefficients in $\mathcal{U}(\mathcal{M}^e).  $
\end{defn}
Note that if $ \mathcal{A} $ is a $ * $-subalgebra of $ \mathcal{U}(\mathcal{M}) $ by  \cite[Corollary 3.7]{t} we have an isomorphism of right $\mathcal{U}(\mathcal{M}^e)  $-modules
$$H^k(\mathcal{A}, \mathcal{U}(\mathcal{M}^e))\cong H_k(\mathcal{A}, \mathcal{U}(\mathcal{M}^e))^\prime\cong (H_k(\mathcal{A},\mathcal{M}^e)\otimes_{\mathcal{M}^e} \mathcal{U}(\mathcal{M}^e))^\prime.$$

Combining this idea with the definition of equivariant Hochschild cohomology due to K. Jensen (\cite[Definition I.1.7]{j}) we arrive at the following (apparently new) notion:
\begin{defn}
	Let $\mathcal{M}$ be a  von Neumann algebra with a finite, faithful,  normal trace $ \tau $ and let $ G \curvearrowright _{\alpha}\mathcal{M} $.
	Suppose that $\tau  $ is invariant under the action $ \alpha $.
	Let $ \mathcal{A} $ be a $ * $-subalgebra of $ \mathcal{M} $ invariant under the action $ \alpha $. We denote the restricted action on $ \mathcal{A} $  again  by $ \alpha $.

	The $ k $-th equivariant $ L^2 $-cohomology of $ \mathcal{A} $ is defined by $$ H_{G}^{(2),k}(\mathcal{A})=H^k_{G}(\mathcal{A}, \mathcal{U}(\mathcal{M}^e)),  $$
	where $H^k_{G}(\mathcal{A}, \mathcal{U}(\mathcal{M}^e))  $ is the $ k $-th  equivariant Hochschild cohomology of $ \mathcal{A} $ with coefficients in $\mathcal{U}(\mathcal{M}^e).  $
	
	Also, the $ k $-th equivariant $ L^2 $-Betti number of $ \mathcal{A} $ is defined by
	$$b_{G}^{(2),k}(\mathcal{A})= \dim _{\mathcal{U}(\mathcal{M}^e)^\beta} H_{G}^{(2),k}(\mathcal{A}).  $$
	Here $ \mathcal{U}(\mathcal{M}^e)^\beta $ denotes the fixed point algebra of $ \mathcal{U}(\mathcal{M}^e) $ with respect to the action $ \beta $ and the dimension function  is the generalized  dimension function described above.

\end{defn}

The algebra $ \mathcal{U}(\mathcal{M}^e) $ above is a right $ \mathcal{U}(\mathcal{M}^e) $-module with its usual multiplication. This gives  $ C^k(\mathcal{A},\mathcal{U}(\mathcal{M}^e) $ a left $ \mathcal{U}(\mathcal{M}^e) $-module structure with the action
$$\xi\cdot T(a_1,... , a_k)=T(a_1,... , a_k)\xi,\quad (\xi\in \mathcal{U}(\mathcal{M}^e), T\in C^k(\mathcal{A},\mathcal{U}(\mathcal{M}^e)).  $$
The restriction of this left action  to $ C^k_G(\mathcal{A},\mathcal{U}(\mathcal{M}^e) $	
is an action of the fixed point algebra $ \mathcal{U}(\mathcal{M}^e)^\beta $ on
$ C^k_G(\mathcal{A},\mathcal{U}(\mathcal{M}^e) $. That is, if $ \xi \in \mathcal{U}(\mathcal{M}^e)^\beta $ and $ T \in  C^k_G(\mathcal{A},\mathcal{U}(\mathcal{M}^e) $, then $ \xi \cdot T \in  C^k_G(\mathcal{A},\mathcal{U}(\mathcal{M}^e). $
The critical fact here is that this action commutes with $ \delta^k_G $ and hence it induces an action of $ \mathcal{U}(\mathcal{M}^e)^\beta  $ on $ H^k_{G}(\mathcal{A}, \mathcal{U}(\mathcal{M}^e) $. This allows us to take the dimension over $ \mathcal{U}(\mathcal{M}^e)^\beta  $ in the above definition.

\subsection{Grading Versus Dual Action}
Graded $ C^* $-algebras and their relation with Fell bundles studied by R. Exel \cite[Definition 3.1]{ex}.

\begin{defn}
	Let $ G $ be a discrete group and let $ \mathcal{A} $ be a $ C^* $-algebra. $ \mathcal{A} $ is said to be $ G $-graded if there is a collection of linearly independent closed subspaces $ (\mathcal{A}_g)_{g \in G} $ of $ \mathcal{A} $ with the following conditions
	
	(i) $\mathcal{A}_g \mathcal{A}_h \subseteq \mathcal{A}_{gh}  $ for all $ g, h \in G $
	
	(ii) $ \mathcal{A}_g^*\subseteq \mathcal{A}_{g^{-1}} $
	
	(iii) $ \mathcal{A} $ is the closure of the direct sum $ \oplus_{g \in G} \mathcal{A}_g $.
	
	This grading is called topological grading if there exists a conditional expectation of
	$ \mathcal{A} $ to $ \mathcal{A}_e  $. Here $ e $ is the identity element of $ G. $
\end{defn}

Similarly one can  give a definition of a graded von Neumann algebras as follows.

Let $ G $ be a discrete group and let $ \mathcal{M} $ be a von Neumann algebra. A $ G $-grading on $ \mathcal{M} $  is a collection of linearly independent closed subspaces $ (\mathcal{M}_g)_{g \in G} $ of $ \mathcal{M} $ with the following conditions

(i) $\mathcal{M}_g \mathcal{M}_h \subseteq \mathcal{M}_{gh}  $ for all $ g, h \in G $

(ii) $ \mathcal{M}_g^*\subseteq \mathcal{M}_{g^{-1}} $

(iii) $ \mathcal{M} $ is the UWOT  closure of the direct sum $ \oplus_{g \in G} \mathcal{M}_g $.

The homology of graded rings is well studied. A good source  is \cite[Section 20]{fe}.  For the rest of this section, we suppose that
$ G $ is a discrete group and $ G^\circ $ its opposite group.

Let $ \mathcal{M} $ be a  von Neumann algebra and
let $ \mathcal{A}=\oplus_{g \in G} \mathcal{A}_g $ be a $ G $-graded *-subalgebra of $ \mathcal{M} $.
This grading of $ \mathcal{A} $ induces a $ (G\otimes G^\circ) $-grading on $ \mathcal{A}\otimes \mathcal{A}^\circ  $ with homogeneous components
$$  (\mathcal{A}\otimes \mathcal{A}^\circ)_g=\{\sum_i a_i\otimes b^\circ_i| a_i \in \mathcal{A}, b^\circ_i\in \mathcal{A}^\circ \:\hbox{and} \;\deg (a_i)+\deg(b^\circ_i)=g  \}. $$
Let $ (P=\{P_{k,*}\},d)\rightarrow (\mathcal{A},0) $ be a graded projective resolution of $ \mathcal{A} $ with the augmentation map $ \sigma $ in the sense of \cite[Section 20]{fe}, that is, a long exact sequence
\begin{equation}
...\xrightarrow{d}P_{1,*}\xrightarrow{d}P_{0,*}\xrightarrow{\sigma}\mathcal{A}\rightarrow 0.
\end{equation}
Then, $ \hbox{Hom}_{\mathcal{A}\otimes A^\circ}(P_{k,*},\mathcal{M}\otimes\mathcal{M}^\circ) $ is a graded module of linear maps and
\begin{equation}\label{6.2}
0\rightarrow \hbox{Hom}_{\mathcal{A}^e}(P_{0,*},\mathcal{M}^e)\xrightarrow{\delta}\hbox{Hom}_{\mathcal{A}^e}(P_{1,*},\mathcal{M}^e)\xrightarrow{\delta} ...
\end{equation}
is a cochain complex of graded modules with the coboundary maps $ \delta $
defined by $ \delta(T)=-(-1)^{\deg (T)} Td$.
We use the standard notation $ \hbox{Hom}_{\mathcal{A}^e}^{k,*}(P,\mathcal{M}^e)  $ instead of $ \hbox{Hom}_{\mathcal{A}^e}(P_{k,*},\mathcal{M}^e) $. Thus the cochain complex (\ref{6.2}) has the form $ (\hbox{Hom}_{\mathcal{A}^e}^{k,*}(P,\mathcal{M}^e), \delta) $.

Now, let us introduce the $ L^2 $-cohomology of graded von Neumann algebras.
\begin{defn}
	For every $ k\geq 0 $, the $ k $-th graded $ L^2 $-cohomology of $\mathcal{A} $ is defined by
	$$H_{grd}^{(2),k}(\mathcal{A}):= H^{k,*}( \hbox{Hom}_{\mathcal{A}^e}(P,\mathcal{M}^e), \delta )=\hbox{Ext}_{\mathcal{A}^e}^{k}(\mathcal{A},\mathcal{M}^e),  $$
	and the $ k $-th graded $ L^2 $-Betti number of $ \mathcal{A} $ is defined by
	$$ \beta_{grd}^k(\mathcal{A})=\dim _{\mathcal{A}^e}H_{grd}^{(2),k}(\mathcal{A}). $$
\end{defn}

An important fact about graded $ C^* $-algebras is the relation between grading on the $ C^* $-algebra and the action induced by this grading on the $ C^* $-algebra.
Raeburn in \cite[Theorem 3]{ra} showed that if $ G $ is an abelian group and if $\mathcal{A}  $ is a topologically $ G $-graded $ C^* $-algebra, then there is a strongly continuous action $ \alpha $ of the dual group $ \widehat{G} $ on $\mathcal{A}  $ such that for each $ a\in \mathcal{A}_g $, $ \alpha_\gamma(a)=\gamma (g)a $ and the conditional expectation $ F :\mathcal{A}\rightarrow \mathcal{A}_e $ is defined by
$ F(a)=\int_{\widehat{G}}\alpha_\gamma (a)d\gamma $ for all $ a\in \mathcal{A} $.

Conversely, let $ \alpha $ be a $ \widehat{G} $-action  on a $ C^* $-algebra $ \mathcal{A}$. For every $ g\in G $ let $ \mathcal{A}_g=\{a\in\mathcal{A}|\alpha_\gamma (a)=\gamma (g)a,\gamma \in  \widehat{G}  \} $ and put $\mathcal{A}_\alpha =\overline{\oplus_{g \in G}\mathcal{A}_g}  $. Then the collection $(\mathcal{A}_g)_{g\in G}  $ is a grading on
the closed subalgebra $ \mathcal{A}_\alpha $ of $ \mathcal{A} $.

\subsection{Concluding Remarks} The idea of equivariant $L^2$-cohomology presented here is at its very beginning and deserves further development. There are a few important features that needs to be further investigated. The first is finding conditions under which the equivariant $L^2$-Betti numbers vanish. Second, it is nice to have a relation between equivariant $L^2$-cohomology and the $L^2$-cohomology of the corresponding crossed-product (or its subalgebras). Third, it is desirable to have a relation between equivariant $L^2$-cohomology and graded $L^2$-cohomology, as actions and gradings are interrelated as shown above. This is not an easy task, as the two notions are defined here via two different approaches. Summing up, the equivariant $L^2$-cohomology theory and its counterparts sound as a promising field of research with many non trivial and challenging problems.


\subsection*{Acknowledgment}
This work was completed with the support of INSF and Tarbiat Modares University.

\end{document}